\documentclass[12pt,reqno]{amsart}
\usepackage{amscd,amsmath,amsopn,amssymb,amsthm,multicol}
\usepackage{hyperref}
\DeclareMathOperator{\Ad}{Ad}

\DeclareMathOperator{\g}{\mathfrak{g}}
\newcommand{\fr}{\mathfrak}
\DeclareMathOperator{\hh}{\mathfrak{h}}

\DeclareMathOperator{\kk}{\mathfrak{k}}
\DeclareMathOperator{\mm}{\mathfrak{m}}

\DeclareMathOperator{\SO}{SO}
\DeclareMathOperator{\s}{S}

\DeclareMathOperator{\Sp}{Sp}
\DeclareMathOperator{\SU}{SU}
\DeclareMathOperator{\U}{U}

\newcommand{\thickline}{\noalign{\hrule height 1pt}}

\newtheorem{theorem}{Theorem}
\newtheorem{lemma}{Lemma}

\newtheorem{proposition}{Proposition}
\newtheorem{defidition}{Definition}
\newtheorem{example}{Example}
\newtheorem{corollary}{Corollary}

\begin{document}

\title[ Homogeneous manifolds whose geodesics are orbits]
{ Homogeneous manifolds whose geodesics are orbits. Recent results and some open problems}
\author{Andreas Arvanitoyeorgos}

\address{University of Patras, Department of Mathematics, GR-26500 Patras, Greece}
\email{arvanito@math.upatras.gr}

\begin{abstract}
 A homogeneous Riemannian manifold $(M=G/K, g)$ is called a space with homogeneous geodesics or a $G$-g.o. space if 
every geodesic  $\gamma (t)$ of $M$ is an orbit of
a one-parameter subgroup of $G$, that is $\gamma(t) = \exp(tX)\cdot o$, for some non zero vector
$X$ in the Lie algebra of $G$. 
We give an exposition on the subject, by presenting techniques that have been used so far and 
a wide selection of previous and recent results.
 We also present some open problems.

 \medskip

\noindent
{\it 2000 Mathematical Subject Classification.} Primary 53C25. Secondary 53C30.

\medskip

\noindent
{\it Key words.} Homogeneous geodesic; g.o. space; invariant metric;  geodesic vector; geodesic lemma; naturally reductive space;  generalized flag manifold; generalized Wallach space; $M$-space; $\delta$-homogeneous space; pseudo-Riemannian manifold; two-step homogeneous geodesic;
 \end{abstract}
\maketitle

\section*{Introduction}
The aim of the present article is to give an exposition on developments about homogeneous geodesics in Riemannian homogeneous spaces, to present various recent results and discuss some open problems.
One of the demanding problems in Riemannian geometry is the description of geodesics.
By making some symmetry assumptions one could expect that certain simplifications may accur.
Let $(M, g)$  be a homogeneous Riemannian manifold, i.e. a connected Riemannian manifold on which the largest connected group $G$ of isometries acts transitively. Then $M$ can be expressed as a homogeneous space $(G/K, g)$, where $K$ is the isotropy group at a fixed point $o$ of $M$.

Motivated by  well known facts such that,  the geodesics in a Lie group $G$ with a bi-invariant metric are the one-parameter subgroups of $G$, or that the  geodesics in a Riemannian symmetric space $G/K$ are orbits of one-parameter subgroups
in $G/K$, it is natural to search for geodesics in a homogeneous space, which are orbits.  More precisely,
 a geodesic $\gamma(t)$ through the origin $o$ of $M=G/K$ is called {\it homogeneous} if it is an orbit of a one-parameter subgroup of $G$, that is
\begin{equation}\label{1}
 \gamma(t)=\mathrm{exp}(tX)\cdot o, \quad t \in \mathbb{R},
 \end{equation}
where $X$ is a non zero vector in the Lie algebra $\mathfrak{g}$ of $G$.  A non zero vector $X \in \mathfrak{g}$ is called a {\it geodesic vector} if the curve (\ref{1}) is a geodesic.
A homogeneous Riemannian manifold $M=G/K$ is called a {\it g.o. space} if all geodesics are homogeneous with respect to the largest connected group of isometries $I_o(M)$.
Since their first systematic study by O. Kowalski and L. Vanhecke in \cite{Ko-Va}, there has been a lot of research related to homogeneous geodesics and g.o spaces and in  various directions. 

\smallskip
Homogeneous geodesics appear quite often in physics as well.
The equation of motion of many systems of classical
mechanics reduces to the geodesic equation in an appropriate Riemannian manifold
$M$. Homogeneous geodesics in $M$ correspond to ``relative equilibriums" of the
corresponding system (cf. \cite{Arn}). 
For further information about relative equilibria in physics we refer to \cite{Ga-Hu-Wi} and references therein.
 In Lorentzian geometry in particular, homogeneous spaces with the property that all their
 {\it null} geodesics are homogeneous, are candidates for constructing solutions to the $11$-dimensional supergravity, which preserve more than $24$ of the available $32$ supersymmetries.
In fact, all Penrose limits, preserving the amount of supersymmetry of such a
solution, must preserve homogeneity.  This is the case for the Penrose limit of a reductive
homogeneous spacetime along a null homogeneous geodesic (\cite{Fi-O'F-Ph-Me}, \cite{Me}, \cite{Phi}).
For a recent mathematical contribution in this topic see \cite{Du4}.

 All naturally reductive spaces  are g.o. spaces (\cite{Ko-No}), but the converse is not true in general. In \cite{Kap} A. Kaplan proved the existence of g.o. spaces that are in no way naturally reductive.  These are generalized Heisenberg groups with two dimensional center.
 Another important class of g.o. spaces are the weakly symmetric spaces.
 These are homogeneous Riemannian manifolds $(M=G/K, g)$  introduced by A. Selberg in \cite{Sel},
with the property that every two points can be interchanged by an isometry of $M$.  
%The classification of weakly symmetric reductive homogeneous spaces was given by O.S. Yakimova in \cite{Ya}.  
In \cite{Be-Ko-Va} J. Berndt, O. Kowalski and L. Vanhecke proved that weakly symmetric spaces are g.o. spaces.
In \cite{Ko-Pr-Va} O. Kowalski, F. Pr\" ufer and L. Vanhecke gave an explicit classification of all naturally reductive spaces up to dimension five. 

The term {\it g.o. space} was introduced by  O. Kowalski and L. Vanhecke in \cite{Ko-Va}, 
where they gave the classification of all g.o. spaces up to dimension six, which are in no way naturally reductive. 
Concerning the existence of homogeneous geodesics in  a homogeneous Riemannian manifold, we recall the following.
 In (\cite{Kaj}) V.V. Kajzer  proved that a Lie group endowed with a left-invariant metric admits at least one homogeneous geodesic. O. Kowalski and J. Szenthe extended this result to all homogeneous Riemannian manifolds (\cite{Ko-Sz}). An extension of this result  to reductive homogeneous  pseudo-Riemannian manifolds was  obtained (\cite{Du-Ko2}, \cite{Phi}).

In \cite{Gor} C. Gordon described g.o. spaces  which are nilmanifolds and in \cite{Tam} H. Tamaru  classified homogeneous g.o. spaces which are fibered over irreducible symmetric spaces.
In  \cite{Du2} and \cite{Du-Ko1} O. Kowalski and Z. Du\v sek investigated  homogeneous geodesics in   Heisenberg groups  and some $H$-type groups. Examples of g.o. spaces in dimension seven were obtained by Du\v sek, O. Kowalski and S. Nik\v cevi\' c in \cite{Du-Ko-Ni}. 

In \cite{Al-Ar} the  author and D.V. Alekseevsky  classified generalized flag manifolds which are   g.o. spaces.
Further, D.V. Alekseevsky and Yu. G. Nikonorov in \cite{Al-Ni} studied the structure of compact g.o. spaces and gave some sufficient conditions for existence and non existence of an invariant metric  with homogeneous geodesics on
a homogeneous space of a compact Lie group. They also gave a classification of compact simply
connected g.o. spaces  of positive Euler characteristic.

 In  \cite{Ko-Ni-Vl} O. Kowalski, S. Nik\v cevi\'c and Z. Vl\' a\v sek studied homogeneous geodesics in homogneous Riemannian manifolds, and  in  \cite{Mar},  \cite{Ca-Ma1} G. Calvaruso  and R. Marinosci studied homogeneous geodesics in three-dimensional Lie groups.
 Homogeneous geodesics were also studied by J. Szenthe in \cite{Sze1}, \cite{Sze2}, \cite{Sze3}, \cite{Sze4}.
 Also, D. Latifi studied homogeneous geodesics in homogeneous Finsler spaces (\cite{Lat}), and the first author  investigated homogeneous  geodesics in  the flag manifold $\SO(2l+1)/\U(l-m)\times \SO(2m+1)$ (\cite{Arv}).

Homogeneous geodesics in the affine setting were studied in \cite{Du2} and \cite{Du-Ko-Vl} 
(and in particular for any non reductive pseudo-Riemannian manifold).

Finally, a class of homogeneous spaces which satisfy the g.o. property are the $\delta$-homogeneous spaces, which were introduced by V. Berestovski\v i and C. Plaut in \cite{Be-Pl}.
These spaces have interesting geometrical properties, but we will not persue here.  We refer to the paper
\cite{Be-Ni} by V. Berestovski\v i and Yu.G. Nikonorov for more information ib this direction.

\smallskip
The paper is organized as follows. In Section 1 we present the basic techniques for finding homogeneous geodesics and detecting if a homogeneous space is a space with homogeneous geodesics (g.o. space).
In Section 2 we present the classification up to dimension $6$ and give examples in dimension $7$.
In Section 3 we discuss homogeneous g.o. spaces which are fibered over irreducible symmetric spaces and in Section 4 we present the classification of generalized flag manifolds which are g.o. spaces.
In Section 5 we present results about another wide class of homogeneous spaces, the generalized Wallach spaces, and in Section 6 we discuss results related to $M$-spaces. These are homogeneous spaces $G/K_1$ so that
$G/(S\times K_1)$ is a generalized flag manifold, where $S$ a torus in a compact simple Lie group $G$.
The pseudo-Riemannian setting is presented in Section 7.  In Section 8 we discuss a generalization of homogeneous geodesics which we call two-step homogeneous geodesics.
These are orbits of the product of two exponential factors.  Finally, in Section 9 we present some open problems.

 \medskip

\noindent
{\bf Acknowledgements.}  The author was supported by Grant \# E.037 from the research committee of the
University of Patras (programme K. Karatheodori). 
He had useful discussions with Zdenek Du\v sek and Giovanni Calvaruso (during his visit in Lecce via the Erasmus+ programme, May 2016).

\section{Homogeneous geodesics in homogeneous Riemannian manifolds - Techniques}\label{Section1}

 Let $(M=G/K, g)$ be a homogeneous space of a compact, connected and semisimple Lie group and let 
 $\mathfrak{g}$,  $\mathfrak{k}$ be the Lie algebras of $G$ and $K$ respecively.
 We consider an orthogonal reductive decomposition
\begin{equation}\label{reductive}
  \mathfrak{g}=\mathfrak{k}\oplus \mathfrak{m}
\end{equation}
with respect to 
 the negative of the Killing form of $\mathfrak{g}$, denoted by $B$. 
 The canonical projection $\pi:G \rightarrow G/K $ induces an isomorphism between the subspace $\mathfrak{m}$ and the tangent space $T_{o}M$ at the identity $o=eK$. 
 A $G$-invariant Riemannian metric $g$   defines a scalar product $\langle \cdot,\cdot \rangle$ on $\mathfrak{m}$ which is $\Ad({K})$-invariant. Then any $\Ad({K})$-invariant scalar product $\langle \cdot,\cdot \rangle$ on $\mathfrak{m}$ can be expressed as $\langle x, y \rangle=B(\Lambda x, y)\; (x, y \in \mathfrak{m})$, where $\Lambda$ is an $\Ad({K})$-equivariant positive definite symmetric operator on $\mathfrak{m}$. Conversely, any such operator $\Lambda$ determines an  $\Ad({K})$-invariant scalar product $\langle x, y \rangle=B(\Lambda x, y)$ on $\mathfrak{m}$, which in turn determines   a $G$-invariant Riemannian metric $g$ on $\mathfrak{m}$. We say that  $\Lambda$ is the {\it operator associated} to the metric $g$, or simply the {\it associated operator}.
Also, a Riemannian metric generated by the scalar product
  product $B$ is called {\it standard metric}.
 
\begin{defidition}\label{go}   A homogeneous Riemannian manifold $(M=G/K, g)$ is called a space with homogeneous geodesics, or $G$-g.o. space if every geodesic $\gamma$ of $M$ is an orbit of a one-parameter subgroup of $G$, that is
  $\gamma (t)=\exp (tX)\cdot o$, for some non zero vector $X\in\mathfrak{g}$.  The invariant metric $g$ is called $G$-g.o. metric.
  If $G$ is the full isometry group, then the $G$-g.o. space is called a manifold with homogeneous geodesics, or a g.o. manifold.
  \end{defidition}

Notice that if all geodesics through the origin $o=eK$ are of the form $\gamma(t)=\exp (tX)\cdot o$, then the geodesics through any other point $a\cdot p$ $(a\in G, p\in M)$ is of the form $a\gamma(t)=\exp(t\Ad(a)X)\cdot (a\cdot p)$.

\begin{defidition} \label {D1}  A non zero vector $X \in \mathfrak{g}$ is called a geodesic vector if the curve (\ref{1}) is a geodesic.
\end{defidition}

All calculations for a g.o. space $G/K$ can be reduced to algebraic computations using geodesic vectors.
These can be computed by using the following fundamental result of the subject, still call it ``lemma" by tradinion:

\begin{lemma}[Geodesic Lemma \cite{Ko-Va}]  \label{L1} A nonzero vector  $X \in \mathfrak{g}$ is a geodesic vector if and only if

\begin{equation}\label{3}
  \langle [X,Y]_{\mathfrak{m}},X_{\mathfrak{m}} \rangle =0,
\end{equation}

for all $Y \in \mathfrak{m}$.
Here the subscript $\mathfrak{m}$ denotes the projection into $\mathfrak{m}$.
\end{lemma}

A useful description of homogeneous geodesics (\ref{1}) is provided by the following :
\begin{proposition}\label{genlemma} \textnormal{(\cite{Al-Ar})}
 Let $(M=G/K, g)$ be a homogeneous Riemannian manifold and $\Lambda$ be the associated operator. Let $a\in \mathfrak{k}$ and $x \in \mathfrak{m}$. Then the following are equivalent:

 (1)\  The orbit $\gamma(t)=\mathrm{exp}t(a+x)\cdot o$ of the one-parameter subgroup $\mathrm{exp}t(a+x)$ through the point $o=eK$ is a geodesic of $M$.

 (2) \ $[a+x, \Lambda x] \in \mathfrak{k}$.

 (3) \ $\langle [a, x], y \rangle = \langle x, [x, y]_{\mathfrak{m}} \rangle \ \mbox {for all} \  y \in \mathfrak{m}$.

 (4) \ $\langle [a+x, y]_{\mathfrak{m}}, x \rangle=0 \  \mbox {for all} \ y \in \mathfrak{m}$.
 \end{proposition}

As a consequence, we obtain the following characterization of g.o. spaces:

 \begin{corollary} [\cite{Al-Ar}]  \label{C1}

 Let $(M=G/K, g)$ be a homogeneous Riemannian manifold. Then $(M=G/K, g)$ is a g.o. space if and only if for every $x\in \mathfrak{m}$  there exists an $a(x) \in \mathfrak{k}$  such that
 $$
 [a(x)+x, \Lambda x] \in \mathfrak{k}.
 $$
 \end{corollary}

Therefore, the property of being a g.o. space $G/K$, depends only on the reductive decomposition  and the 
$G$-invariant metric metric $g$ on $\mathfrak{m}$. 
That is, if $(M = G/H, g)$ is a g.o. space,
then any locally isomorphic homogeneous Riemannian space $(M = G/H, g')$
is a g.o. space.
Also, a direct product of Riemannian manifolds is a manifold with homogeneous geodesics if and
only if each factor is a manifold with homogeneous geodesics.

\medskip
In order to find all homogeneous geodesics in a
  homogeneous Riemannian manifold $(M=G/K, g)$ it suffices
  to   find a decomposition of the form (\ref{reductive}) and look for geodesic vectors of  the form
  \begin{equation}\label{4}
  X=\sum^{s}_{i=1}x_{i}e_{i}+\sum^{l}_{j=1}a_{j}A_{j}.
  \end{equation}
  Here $\{e_{i}: i=1,2, \dots, s\}$ is a convenient basis of $\mathfrak{m}$ and $\{A_{j}: j=1,2,\dots,l\}$ is a  basis of $\mathfrak{k}$.
  By substituting $X=e_i$ $(i=1, \dots ,s)$ into equation (\ref{3}) we obtain a system of linear algebraic equations for the variables $x_i$ and $a_j$. 
The geodesic vectors correspond to those solutions for which $x_1, \dots , x_s$ are not all equal to zero.
 For some applications of this method we refer to \cite{Ko-Ni-Vl} and \cite{Mar}.
  Also, $(M=G/K, g)$ is a g.o. space if and only if for every non zero $s$-tuple $(x_1, \dots, x_s)$ there is an
  $l$-tuple $(a_1, \dots, a_l)$ satisfying all quadratic equations. 
 
 A useful technique used for the characterization of Riemannian g.o. spaces is based on the concept of the geodesic graph, originally introduced in \cite{Sze1}.  We first need the following definition.
 
 \begin{defidition}\label{nr1}  A Riemannian homogeneous space $(G/K, g)$ is called naturally reductive if there exists a reductive decomposition (\ref{reductive}) of $\g$ such that
 \begin{equation}\label{nr}
 \langle [X, Z]_{\mm}, Y\rangle+\langle X, [Z, Y]_{\mm}\rangle =0, \quad\mbox{for all}\ \ X, Y, Z\in \mm.
 \end{equation}
 \end{defidition} 

It is well known that condition (\ref{nr}) implies that all geodesics in $G/K$ are homogeneous (e.g. \cite{ON}).
\medskip

\begin{defidition}\label{nr2} A homogeneous Riemannian manifold $(M, g)$ is naturally reductive if there exists a transitive group $G$ of isometries for which the correseponding Riemannian homogeneous space $(G/K, g)$ is naturally reductive in the sense of Definition \ref{nr1}.
\end{defidition}
Therefore, it could be possible that a homogeneous space $M=G/K$ is not naturally reductive for some group $G\in I_0(M)$ (the connected component of the full isometry group of $M$), but it is naturally reductive
if we write $M=G'/K'$ for some larger group of isometries $G'\subset I_0(M)$.

Let $(M=G/K,g)$ be a g.o. space and let $\g = \kk\oplus\mm$ be an $\Ad(K)$-invariant decomposition.  Then

(1) There exists an $\Ad(H)$-equivariant map $\eta:\mm\to\kk$ (a {\it geodesic graph}) such that for any $X\in\mm\setminus\{0\}$, the curve $\exp t(X+\eta(X))\cdot o$ is a geodesic.

(2) A geodesic graph is either linear (which is equivalent to natural reductivity with respect to some
$\Ad(K)$-invariant decomposition $\g=\kk\oplus\mm '$) or it is non differentiable at the origin $o$.

It can be shown (\cite{Ko-Ni}) that a geodesic graph (for  a g.o. space) is uniquely determined by fixing an $\Ad(H)$-invariant scalar product on $\kk$.
Examples of g.o. spaces by using geodesic graphs are given in \cite{Du5}, \cite{Du-Ko-Ni}, and \cite{Ko-Ni}.
Conversely, the property (1) implies that $G/K$ is a g.o. space.

Another technique for producing g.o. metrics was given by C. Gordon as shown below:

\begin{proposition}\textnormal{(\cite{Gor}, \cite{Tam})}\label{P5}
Let $G$ be a connected semisimple Lie group and $H\supset K$ be  compact Lie subgroups in $G$.  Let  $M_{F}$ and $M_{C}$ be the  tangent spaces of $F=H/K$ and  $C=G/H$ respectively. Then the metric $g_{a, b}= a B\mid_{M_{F}}+b B\mid_{M_{C}}, (a, b \in \mathbb{R}^+)$ is a g.o. metric on $G/K$ if and only if  for any $v_F\in M_{F}$, $v_C \in M_{C}$ there exists $X\in \mathfrak{k}$ such that
\begin{equation*}
[X, v_F]=[X+v_F, v_C]=0.
\end{equation*}
\end{proposition} 

Actually, Gordon proved a more general result based on description of naturally reductive left-invariant metrics on compact Lie groups given by J.E. D'Atri and W. Ziller in \cite{D'A-Zi}.

\section{Low dimensional examples}

The problem of a complete classification of g.o. manifolds is open. Even the classification all g.o. metrics on a given Riemannian homogeneous space is not trivial (cf. for example \cite{Nik2}).
A complete classification is known up to dimesion $6$, given by O. Kowalski and L. Vanhecke: 

\begin{theorem}\textnormal{(\cite{Ko-Va})}
\textnormal{1)} All Riemannian g.o. spaces of dimension up to $4$ are naturally reductive.

\textnormal{2)} Every $5$-dimensional Riemannian g.o. space is either naturally reductive, or of isotropy type $\SU(2)$.

\textnormal{3)} Every $6$-dimensional Riemannian g.o. space is either naturally reductive or one of the following:

\  \textnormal{(a)}  A two-step nilpotent Lie group with two-dimensional center, equipped with a left-invariant
Riemannian metric such that the maximal connected isotropy group is isomorphic to either  $\SU(2)$ or 
$\U(2)$. 
Corresponding g.o. metrics depend on three real parameters.

\ \textnormal{(b)} The universal covering space of a homogeneous Riemannian manifold of the form
$(M = \SO(5)/\U(2), g)$ or $(M = \SO(4, 1)/\U(2), g)$, where $\SO(5)$ or $\SO(4, 1)$ is the identity component
of the full isometry group, respectively. In each case, all corresponding invariant metrics g.o. metrics $g$ depend on two real parameters.
\end{theorem}
As pointed out by the authors in \cite[p. 190]{Ko-Va}, the g.o. spaces (a) and (b) are {\it in no way naturally reductive} in the following sence: whatever the representation of $(M,g)$ as a quotient of the form $G'/K'$, where $G'$ is a connected transitive group of isometries of $(M,g)$, and whatever is the $\Ad(K)$-invariant decomposition $\g '=\kk '\oplus\mm '$, the curve $\gamma(t)=\exp (tX)\cdot o$ is never a geodesic (for any $X\in\mm\setminus\{0\}$). 

\medskip
The first $7$-dimensional example of a g.o. manifold was given by C. Gordon in \cite{Gor}.  This is a nilmanifold (i.e. a connected Riemannian manifold admitting a transitive nilpotent group of isometries), and it was obtained under a general construction of g.o. nilmanifolds.
It took some time until some more $7$-dimensional examples were given.  In \cite{Du-Ko-Ni} 
Z. Du\v sek, O. Kowalski and S. Nik\v cevi\' c gave families of $7$-dimensional g.o. metrics.  Their main result is the following:

\begin{theorem}\textnormal{(\cite{Du-Ko-Ni})} 
On the $7$-dimensional homogeneous space $G/K = (\SO(5) \times \SO(2))/\U(2)$ (or $G/H =
(\SO(4, 1) \times \SO(2))/\U(2)$) there is a family ${g_{p,q}}$ of invariant metrics depending on two
parameters $p, q$ (where the pairs $(p, q)$ fill in an open subset of the plane) such that each homogeneous
Riemannian manifold $(G/H, g_{p,q})$ is a locally irreducible and not naturally reductive Riemannian g.o.
manifold.
\end{theorem}

\section{Fibrations over symmetric spaces}

In the work \cite{Tam} H. Tamaru classified homogeneous spaces $M=G/K$ satisfying the following properties:
(i) $M$ is fibered over irreducuble symmetric spaces $G/H$ and (ii) certain $G$-invariant metrics on $M$ are $G$-g.o. metrics.
More precisely, for $G$ connected and semisimple, and $H, K$ compact with $G\supset H\supset K$, he considered the fibration
$$
F=H/K\to M=G/K\to B=G/H
$$
 and the $G$-invariant metrics $g_{a, b}$ on $M$ determined by the scalar products
 $$
 \langle\ ,\ \rangle =\left. aB\right|_{\mathfrak{f}}+\left. bB\right|_{\mathfrak{b}}, \quad a,b>0.
 $$
 Here $\mathfrak{f}$ and $\mathfrak{b}$ are the tangent spaces of $F$ and $B$ respectively, so that
 the tangent space of $M$ at the origin is identified with $\mathfrak{f}\oplus\mathfrak{b}$.
 By using results from polar representations, he classifed all  triplets $(G, H, K)$ so that the metrics $g_{a, b}$ are $G$-g.o. metrics.
 The  triplets of Lie algebras $(\g, \hh, \kk)$ 
so that $(\g, \hh)$ is a symmetric pair and $(\g, \kk)$ corresponds to a $G$-g.o. space $G/K$,
 are listed in Table \ref{tamaru}.

\begin{table}  
 \begin{tabular}{|c||c|c|c|}
  \hline
  & $\g$ & $\hh$ & $\kk$ \\
  \hline
  1& $\mathfrak{so}(2n+1), n\ge 2$ & $\mathfrak{so}(2n)$ & $\mathfrak{u}(n)$\\
 \hline
  2& $\mathfrak{so}(4n+1), n\ge 1$ & $\mathfrak{so}(4n)$ & $\mathfrak{su}(2n)$\\
 \hline
  3& $\mathfrak{so}(8)$ & $\mathfrak{so}(7)$ & $\mathfrak{g}_2$\\
 \hline
  4& $\mathfrak{so}(9)$ & $\mathfrak{so}(8)$ & $\mathfrak{so}(7)$\\
 \hline
  5&$\mathfrak{su}(n+1), n\ge 2$ & $\mathfrak{u}(n)$ & $\mathfrak{su}(n)$\\
 \hline
  6&$\mathfrak{su}(2n+1), n\ge 2$ & $\mathfrak{u}(2n)$ & $\mathfrak{u}(1)\oplus\mathfrak{sp}(n)$\\
 \hline
  7&$\mathfrak{su}(2n+1), n\ge 2$ & $\mathfrak{u}(2n)$ & $\mathfrak{sp}(n)$\\
 \hline
  8&$\mathfrak{sp}(n+1), n\ge 1$ & $\mathfrak{sp}(1)\oplus\mathfrak{sp}(n)$ & $\mathfrak{u}(1)\oplus\mathfrak{sp}(n)$\\
 \hline
  9&$\mathfrak{sp}(n+1), n\ge 1$ & $\mathfrak{sp}(1)\oplus\mathfrak{sp}(n)$ & $\mathfrak{sp}(n)$\\
 \hline
  10&$\mathfrak{su}(2r+n), r\ge 2, n\ge 1$ & $\mathfrak{su}(r)\oplus\mathfrak{su}(r+n)\oplus\mathbb{R}$ & 
    $\mathfrak{su}(r)\oplus\mathfrak{su}(r+n)$\\
    \hline
 11& $\mathfrak{so}(4n+2), n\ge 2$ & $\mathfrak{u}(2n+1)$ & $\mathfrak{su}(2n+1)$\\
 \hline
 12& $\mathfrak{e}_6$ & $\mathbb{R}\oplus\mathfrak{so}(10)$ & $\mathfrak{so}(10)$\\
 \hline
 13& $\mathfrak{so}(9)$ & $\mathfrak{so}(7)\oplus\mathfrak{so}(2)$ & $\mathfrak{g}_2\oplus\mathfrak{so}(2)$\\
 \hline
  14&$\mathfrak{so}(10)$ & $\mathfrak{so}(8)\oplus\mathfrak{so}(2)$ & $spin(7)\oplus\mathfrak{so}(2)$\\
 \hline
  15&$\mathfrak{so}(11)$ & $\mathfrak{so}(8)\oplus\mathfrak{so}(3)$ & $spin(7)\oplus\mathfrak{so}(3)$\\
 \hline
 \end{tabular}
 \medskip
 \caption{Riemannian g.o. spaces $G/K$ fibered over irreducible symmetric spaces $G/H$ (\cite{Tam}).}\label{tamaru}
 \end{table}
 
 %Furthermore, as shown in \cite{Tam} the associated homogeneous spaces corresponding to the pairs
 %6 and 7 in Table \ref{tamaru} are weakly symmetric spaces with respect to any $\SU(2n+1)$-invariant metric.
 
 \section{Generalized flag manifolds}
 
 In the work \cite{Al-Ar} D.V. Alekseevsky and the author 
 classified generalized flag manifolds with homogeneous geodesics.
   Recall that a generalized flag manifold is a homogeneous space
    $G/K$ which is  an adjoint orbit of a  compact semisimple Lie group $G$.
    Equivalently,  the isotropy subgroup  $K$ is the centralizer of a torus
    (i.e. a maximal abelian subgroup) in $G$.
    We assume that $G$  acts effectively on $M$.
    A flag manifold $M=G/K$ is simply connected and has
    the
 canonically defined decomposition
 $  M = G/K = G_1/K_1 \times G_2/K_2 \times \cdots \times G_n/K_n $,
 where  $G_1, \dots , G_n$ are simple factors of the (connected) Lie group $G$.
  This decomposition is the de Rham decomposition  of $M$ equipped with a
 $G$-invariant metric $g$. In particular, $(M,g)$ is a g.o. space
 if
 and only if each factor $(M_i = G_i/K_i, g_i = \left.g \right |_{M_i})$
 is a g.o. space.
 This reduces the problem of the description of $G$-invariant metrics with
 homogeneous geodesics in a flag manifold $M=G/K$ to the case when the group
 $G$ is simple.
 
 Flag manifolds $M=G/K$ with $G$ a simple Lie group can be classified in terms of their
 {\it painted Dynkin diagrams}.  It turns out that
 for each classical Lie group there is an infinite series of flag manifolds,
 and for each of the exceptional Lie groups $G_2$, $F_4$, $E_6$, $E_7$, and $E_8$
  there are 3, 11, 16, 31, and 40 non equivalent flag manifolds respectively (eg. \cite{Al}, \cite{Bor}). 
 An important invariant of  flag manifolds is their set of $T$-roots $R_T$.
 This is defined as the restriction of the
 root system $R$ of $\g$ to the center $\mathfrak{t}$ of the stability subalgebra $\kk$ of
 $K$.  
 In \cite{Al-Ar} we
	defined the notion of 
	{\it connected component}
	of $R_T$, namely two $T$-roots are in the same component if they can be
	connected by a chain of $T$-roots whose sum or difference is also a $T$-root.
	The set $R_T$ is called {\it connected} if it
	has only one connected component.
	
	\begin{theorem}\textnormal{(\cite{Al-Ar})} If the set of $T$-roots is connected then the standard metric on 
	 $M=G/K$ is the only $G$-invariant metric (up to scalar) which is a g.o. metric.
	 \end{theorem}
	
	Hence, for a flag manifold $M=G/K$ ($G$ simple), a $G$-invariant g.o. metric 
may exist, only when $R_T$ is not connected,  
	so we only need to study those
  flag manifolds.
  It turns out  that the system of $T$-roots is not connected only for
 three
 infinite series of a classical Lie group (namely the spaces
 $\SO(2\ell +1)/\U(\ell-m)\cdot \SO(2m+1)$, $\Sp(\ell)/\U(\ell -m)\cdot \Sp(m)$, and
 $\SO(2\ell)/\U(\ell-m)\cdot \SO(2m)$),
 and for 10 flag manifolds of an exceptional Lie group.
 An perpective of the above theorem is given by the following theorem:
 
 \begin{theorem}\textnormal{(\cite{Al-Ar})}\label{AA2} Let $M=G/K$ be a flag manifold of a simple Lie group.
 Then the set of $T$-roots is not connected if and only if the isotropy 
 representation of $M$ consists of two irreducible (non-equivalent) components.
 \end{theorem}
 
 Therefore,
 the problem of the description of $G$-invariant metrics on flag manifolds with homogeneous
 geodesics reduces substantially to the study of this short list of
 prospective flag manifolds.  
 To this end, we used the classification Table \ref{tamaru} of the work of H. Tamaru (\cite{Tam}). 
 Since any flag manifold can be fibered over a symmetric space (\cite{Bur}),
  then by using Theorem \ref{AA2}
 we obtain that the only flag manifolds which are in Table \ref{tamaru} are
 $\SO(2\ell +1)/\U(\ell)$ and $\Sp(\ell)/\U(1)\cdot \Sp(\ell -1)$.

 On the other hand, in \cite{Ak-Vin} D.N. Akhiezer and E.B. Vinberg classified all compact
weakly symmetric spaces.  Their classification shows that the only flag
manifolds which are weakly symmetric spaces are
$\SO(2\ell +1)/\U(\ell)$ and $C(1, \ell -1)=\Sp(\ell)/U(1)\cdot
\Sp(\ell -1)$.
This implies that any $\SO(2\ell +1)$-invariant metric $g_\lambda$ on
$\SO(2\ell +1)/\U(\ell)$ (depending, up to scale, on one real parameter $\lambda$)
is weakly symmetric, hence it has homogeneous geodesics.
Similarly for any $\Sp(\ell)$-invariant metric $g_\lambda$ on $\Sp(\ell)/U(1)\cdot
\Sp(\ell -1)$.
In fact, the action of the group $\SO(2\ell +1)$ on $\SO(2\ell +1)/\U(\ell)$ can be
extented to the action of the group $\SO(2\ell +2)$ with isotropy subgroup
$U(\ell +1)$, which preserves the complex structure and the standard
invariant metric $g_0$ (which corresponds to $\lambda =1$).
Hence, the Riemannian flag manifold $(\SO(2\ell +1)/\U(\ell), g_0)$ is isometric
to the Hermitian symmetric space $\mbox{Com}(\mathbb{R} ^{2\ell +2})=\SO(2\ell +2)/\U(\ell +1)$
of all complex structures in $\mathbb{R} ^{2\ell +2}$.
Similarly, the action of the group $\Sp(\ell)$ on $\Sp(\ell)/U(1)\cdot
\Sp(\ell -1)$
can be extended
to the action of the group $\SU(2\ell)$ with isotropy subgroup
$\s(\U(1)\times \U(2\ell -1))$. 
As a consequence of the above we obtain the following:

\begin{theorem}\textnormal{(\cite{Al-Ar})}\label{AA3}  The only flag manifolds $M=G/K$ of a simple Lie group
$G$ admiting a non naturally reductive $G$-invariant metric $g$ with
homogeneous geodesics are the manifolds
$\SO(2\ell +1)/U(\ell)$ and
$\Sp(\ell)/\U(1)\cdot \Sp(\ell -1)$ $(\ell\ge 2)$, which admit (up to scale)
a one-parameter family of $\SO(2\ell +1)$ (resp. $\Sp(\ell)$)-invariant metrics $g_\lambda$.
Moreover, these manifolds are weakly symmetric spaces
for $\lambda >0$, and they are symmetric spaces with respect to
$\mbox{Isom}(g_\lambda)$ if and only if $\lambda =1$, i.e.
$g_\lambda$ is a multiple of the standard metric.
\end{theorem}

Note that for $\ell =2$ we obtain $\Sp(2)/U(1)\cdot \Sp(1)\cong \SO(5)/\U(2)$,
where the second quotient 
is an example of g.o. space in \cite{Ko-Va}
which is not naturally reductive.

Finally, we mention a remarkable coincidence
between Theorem \ref{AA3} and a resut  
by F. Podest\` a and G. Thorbergsson in \cite{Po-Th}, where they studied
coisotropic actions on flag manifolds.
One of their theorems states that 
if $M=G/K$ is a flag manifold of a simple Lie group then the action of $K$ on $M$
is coisotropic, if and only if $M$ is up to local isomorphy either a Hermitian
symmetric space, or one of the spaces obtained in Theorem \ref{AA3}.

 \section{Generalized Wallach spaces}

Let $G/K$ be a compact homogeneous space with connected compact semisimple Lie group $G$ and a compact subgroup $K$ with reductive decomposition $\g =\kk\oplus\mm$.
  Then $G/K$ is called {\it generalized Wallach space} 
(known before as three-locally-symmetric spaces, cf. \cite{Lo-Ni-Fi})  
  if the  module $\mathfrak{m}$  decomposes into a direct sum of three $\Ad(K)$-invariant irreducible modules pairwise orthogonal with respect to $B$, i.e.
  $
  \mathfrak{m}=\mathfrak{m}_1\oplus\mathfrak{m}_2\oplus\mathfrak{m}_3,
  $
  such that
 $
  [\mathfrak{m}_i, \mathfrak{m}_i]\subset\mathfrak{k} \quad  i=1,2,3$.
   Every generalized Wallach space admits a three parameter family of invariant Riemannian metrics determined by $\Ad({K})$-invariant inner products
 $
  \langle \cdot,\cdot\rangle=\lambda_1B(\cdot,\cdot)\mid_{\mathfrak{m}_1}+\lambda_2B(\cdot,\cdot)\mid_{\mathfrak{m}_2}+\lambda_3B(\cdot,\cdot )\mid_{\mathfrak{m}_3}$,
  where $\lambda_1, \lambda_2, \lambda_3$ are positive real numbers.
 
  A classification of generalized Wallach spaces  was recently obtained by Yu.G. Nikoronov (\cite{Nik3}) and  Z. Chen, Y. Kang, K. Liang (\cite{Ch-Ka-Li}) as follows:

\begin{theorem} [\cite{Nik3}, \cite{Ch-Ka-Li}]\label{T1} Let $G/K$ be a connected and simply connected compact homogeneous space. Then $G/K$ is a generalized Wallach space if and only if it is one of the following types:

  1) $G/K$  is a direct product of three irreducible symmetric spaces of compact type.

  2) The group is simple and the pair $(\mathfrak{g},\mathfrak{k})$ is   one of the pairs in Table 1.

  3) $G=F\times F\times F\times F$ and $H=diag(F)\subset G$ for some connected, compact, simple Lie group $F$, with the following description on the Lie algebra level:
  $$
  (\mathfrak{g},\mathfrak{k})=(\mathfrak{f}\oplus \mathfrak{f}\oplus\mathfrak{f}\oplus \mathfrak{f}, {\rm diag}(\mathfrak{f}))=\{(X,X,X,X)\mid X \in f\},
  $$
  where $\mathfrak{f}$ is the Lie algebra of $F$, and (up to permutation) $\mathfrak{m}_1=\{(X,X,-X,-X)\mid X \in f\}$, $\mathfrak{m}_2=\{(X,-X,X,-X)\mid X \in f\}$,  $\mathfrak{m}_3=\{(X,-X,-X,X)\mid X \in f\}$.
  \end{theorem}

\medskip
\begin{center}
\begin{tabular}{|c|c|c|c|}
 \hline
         $\mathfrak{g}$  & $\mathfrak{k}$   & $\mathfrak{g}$  & $\mathfrak{k}$  \\
     \thickline
        $\mathfrak{so}(k+l+m)$ & $\mathfrak{so}(k)\oplus \mathfrak{so}(l)\oplus \mathfrak{so}(m)$ & $\mathfrak{e}_7$&$\mathfrak{so}(8)\oplus 3 \mathfrak{sp}(1)$   \\
    \thickline
       $\mathfrak{su}(k+l+m)$ & $\mathfrak{su}(k)\oplus su(l)\oplus \mathfrak{su}(m)$  &$\mathfrak{e}_7$ &$\mathfrak{su}(6)\oplus \mathfrak{sp}(1)\oplus \mathbb{R}$ \\
    \thickline
     $\mathfrak{sp}(k+l+m)$  &$\mathfrak{sp}(k)\oplus \mathfrak{sp}(l)\oplus \mathfrak{sp}(m)$  &$\mathfrak{e}_7$&$\mathfrak{so}(8)$ \\
   \thickline
        $\mathfrak{su}(2l), l\geq 2$  & $\mathfrak{u}(l)$ & $\mathfrak{e}_8$  &$\mathfrak{so}(12)\oplus 2 \mathfrak{sp}(1)$ \\
    \thickline
       $\mathfrak{so}(2l), l\geq 4$ &$\mathfrak{u}(l)\oplus \mathfrak{u}(l-1)$ &$\mathfrak{e}_8$ & $\mathfrak{so}(8)\oplus \mathfrak{so}(8)$ \\
    \thickline
       $\mathfrak{e}_6$& $\mathfrak{su}(4)\oplus 2 \mathfrak{sp}(1)\oplus \mathbb{R}$ &$\mathfrak{f}_4$&$\mathfrak{so}(5)\oplus 2 \mathfrak{sp}(1)$ \\
    \thickline
        $\mathfrak{e}_6$ & $\mathfrak{so}(8)\oplus \mathbb{R}^2$ &$\mathfrak{f}_4$ &$\mathfrak{so}(8)$ \\
     \thickline
        $\mathfrak{e}_6$&  $\mathfrak{sp}(3)\oplus \mathfrak{sp}(1)$&&\\
  \hline
 \end{tabular}
 \end{center}
 \smallskip
 \begin{center}
  Table 2. \ {\small The pairs $(\mathfrak{g}, \mathfrak{k})$ corresponding to generalized Wallach spaces $G/K$ with $G$ simple
  (\cite{Nik3}).}
\end{center}

In \cite{Arv-Wa1} Yu Wang and the author investigated which of the families of spaces listed in Theorem \ref{T1} are g.o. spaces.
By applying the method of searching for geodesics vectors shown at the end of  Section \ref{Section1} we obtained the following:

\begin{theorem}\label{T3}\textnormal{\cite{Arv-Wa1}} Let $(G/K, g)$ be a generalized Wallach space as listed in Theorem \ref{T1}. Then

1)\ If $(G/K, g)$ is a space  of  type 1)  then this is a  g.o. space  for any $\mathrm{Ad}(K)$-invariant Riemannian metric.

2) \ If  $(G/K, g)$ is a space  of  type 2) or 3) then this  is a g.o. space if and only if $g$ is the standard metric.
\end{theorem}

 However, in order to find all  homogeneous geodesics in $G/K$ it  suffices  to find all the real solutions of  a system  of $\dim\mathfrak{m}_1+\dim\mathfrak{m}_2+\dim\mathfrak{m}_3$ quadratic equations.

 By Theorem \ref{T3} we only need to consider homogeneous geodesics for spaces of types 2) and 3) given in Theorem \ref{T1}, for the metric $(\lambda_{1}, \lambda_{2}, \lambda_{3})$, where  at least two of  $\lambda_{1}, \lambda_{2}, \lambda_{3}$ are  different.  This is not easy in general. 
We obtained all homogeneous geodesics (for various values of the parameters $\lambda_1, \lambda _2, \lambda _3$ for
 the generalized Wallach space $SU(2)/\{e\}$,
  hence recovering a result on R.A. Marinosci (\cite[p. 266]{Mar}),
and for the Stiefel manifolds
 $\SO(n)/\SO(n-2)$, ($n\ge 4$).
 
 \section{$M$-spaces}
 
  Let $G/K$ be a generalized flag manifold with
$K=C(S)=S\times K_1$,
where $S$ is a torus in a compact simple Lie group $G$ and $K_1$ is the semisimple part of $K$.
Then the {\it associated $M$-space} is the homogeneous space $G/K_1$.  These spaces were introduced and studied by H.C. Wang in \cite{Wan}.

 In the work \cite{Arv-Wa-Zh} Y. Wang, G. Zhao and the author 
  investigated homogeneous geodesics in a class of homogeneous spaces called $M$-spaces. 
We proved that for various classes of $M$-spaces, the only g.o. metric is the standard metric.
For other classes of $M$-spaces we give either necessary or necessary and sufficient conditions so that a $G$-invariant metric on $G/K_1$ is a g.o. metric.
The analysis is based on properties of the isotropy representation
$\mathfrak{m}=\mathfrak{m}_1\oplus \cdots\oplus \mathfrak{m}_s$ of the flag manifold $G/K$, in particular on the dimension of the submodules $\mathfrak{m}_i$.
We summarize these results below.

 Let $\mathfrak{g}$  and  $\mathfrak{k}$ be the Lie algebras of the Lie groups $G$ and $K$ respectively. Let $\mathfrak{g}=\mathfrak{k}\oplus \mathfrak{m}$ be an $\Ad(K)$-invariant reductive decomposition of the Lie algebra $\mathfrak{g}$, where $\mathfrak{m}\cong T_o(G/K)$.  This is orthogonal with respect to $B=-$Killing from on $\mathfrak{g}$.
Assume that
 \begin{equation}\label{2}
 \mathfrak{m}=\mathfrak{m}_1\oplus \cdots\oplus \mathfrak{m}_s
 \end{equation}
  is a $B$-orthogonal decomposition of $\mathfrak{m}$ into pairwise inequivalent irreducible $\mathrm{ad}(\mathfrak{k})$-modules.

  Let $G/K_1$ be the corresponding $M$-space  and
   $\mathfrak{s}$ and $ \mathfrak{k}_1$ be  the Lie algebras of $S$ and $K_1$ respectively. We denote by $\mathfrak{n}$  the tangent space $T_o(G/K_1)$, where $o=eK_1$.  A $G$-invariant metric $g$  on $G/K_1$ induces a scalar product $\langle \cdot,\cdot \rangle$ on $\mathfrak{n}$ which is $\Ad({K}_1)$-invariant.  Such an  $\Ad({K}_1)$-invariant scalar product $\langle \cdot,\cdot \rangle$ on $\mathfrak{n}$ can be expressed as $\langle x, y \rangle=B(\Lambda x, y)\; (x, y \in \mathfrak{n})$, where $\Lambda$ is the $\Ad({K}_1)$-equivariant positive definite symmetric operator on $\mathfrak{n}$.

\smallskip
The main results are the following:

\begin{theorem} \label{T1}\textnormal{(\cite{Arv-Wa-Zh})} Let  $G/K$ be a generalized flag manifold  with $s\geq 3 $ in the decomposition  (\ref{reductive}). Let $G/K_1$ be the corresponding $M$-space.

 1) If  $\dim \mathfrak{m}_i\neq 2 \ (i=1,\dots, s)$ and   $(G/K_1, g)$  is a  g.o. space,  then

 $$g=\langle \cdot,\cdot\rangle=\mu B(\cdot,\cdot)\mid_{\mathfrak{s}}+\lambda B(\cdot,\cdot)\mid_{\mathfrak{m}_1\oplus\mathfrak{m}_2\oplus \cdots \oplus \mathfrak{m}_s}, \ (\mu, \lambda >0).
 $$

 2) If there exists  some $j\in  \{1,\dots, s\}$  such that $\dim \mathfrak{m}_{j}= 2$, then $(G/K_1, g)$  is a  g.o. space if and only if $g$ is the standard metric.
\end{theorem}

%\begin{corollary}\label{P4}\label{C1}\textnormal{(\cite{Arv-Wa-Zh})} 
%Let $(G/K_1, g)$ be a $M$-space, where $g=\mu B(\cdot,\cdot)\mid_{\mathfrak{s}}+\lambda B(\cdot,\cdot)\mid_{\mathfrak{m}}$. Then  $(G/K_1, g)$ is a g.o. space for every $\mu, \lambda > 0$ if and only if for every $V\in \mathfrak{s}$  and $X\in \mathfrak{m}$ there exists  $k\in \mathfrak{k}_1$ such that
%\begin{equation*}\label{23}
%[k+V,X]=0.
%\end{equation*}
%\end{corollary}

 \begin{theorem}\label{T2} \textnormal{(\cite{Arv-Wa-Zh})}
 Let $G/K$ be a generalized flag manifold with two isotropy summands $\mathfrak{m}=\mathfrak{m}_1\oplus\mathfrak{m}_2$ and let $(G/K_1, g)$ be the corresponding $M$-space.

1) If $\dim \mathfrak{m}_2=2$ then the standard metric is the only g.o. metric on $M$, unless $G/K_1=\SO(5)/\SU(2)$ or $\Sp(n)/\Sp(n-1)$.
 
2) If  $\dim \mathfrak{m}_i\ne 2$ ($i=1,2$) and the corresponding $M$-space
$(G/K_1, g)$ is a g.o. space, then $g=\langle\cdot ,\cdot\rangle =
\mu B(\cdot,\cdot)\mid_{\mathfrak{s}}+\lambda B(\cdot,\cdot)\mid_{\mathfrak{m}_1\oplus\mathfrak{m}_2}$, 
$(\mu, \lambda >0)$, unless $G/K_1=\SO(2n+1)/\SU(n)$, $n>2$.
 \end{theorem}
However, the spaces $\SO(5)/\SU(2)$ and $\Sp(n)/\Sp(n-1)$ are
included in Tamaru's Table \ref{tamaru}, therefore they admit g.o. metrics.

\section{Homogeneous geodesics in pseudo-Riemannian manifolds}
It is well known that any homogeneous Riemannian manifold
is reductive, but this is not the case for pseudo-Riemannian manifolds in general.
In fact, there exist homogeneous pseudo-Riemannian manifolds which do
not admit any reductive decomposition.  Therefore, there is a dichotomy in the study of geometrical problems
between reductive and non reductive pseudo-Riemannian manifolds.
Due to the existence of null vectors in a pseudo-Riemannian manifold the definition of a homogeneous
 geodesic $\gamma(t)=\exp (tX)\cdot o$ needs to be modified by requiring that $\nabla _{\dot{\gamma}}\dot{\gamma}  =k(\gamma)\dot{\gamma}$ (see also relevant discussion in \cite[pp. 90-91]{Mi}).
It turns out that $k(\gamma)$ is a constant function (cf. \cite{Du-Ko2}.

Even though an algebraic characterization of geodesic vectors (that is an analogue of the geodesic Lemma \ref{L1}) was known to physicists (\cite{Fi-O'F-Ph-Me}, \cite{Phi}), a formal proof was  given  by Z. Du\v sek and O. Kowalski in \cite{Du-Ko2}.  

\begin{lemma}[\cite{Du-Ko2}]  \label{L2}
Let $M = G/H$  be a reductive homogeneous
pseudo-Riemannian space with  reductive decomposition $\fr{g} = \fr{m}\oplus \fr{h}$, and $X\in\fr{g}$.  Then the curve
$\gamma (t)=\exp(tX)\cdot o$ is a 
geodesic curve with respect to some parameter $s$ if and only if
$$
\langle [V, Z]_\fr{m}, V_\fr{m}\rangle = k\langle V_\fr{m}, Z_\fr{m}\rangle,\ \ 
\mbox{for all}\ \  Z \in\fr{m},
$$
where $k$ is some real constant.
Moreover, if $k = 0$, then $t$ is an affine parameter for this geodesic. 
If $k\ne  0$,
then $s = e^{kt}$ is an affine parameter for the geodesic. This occurs
only if the curve $\gamma (t)$ is a null curve in a (properly) pseudo-Riemannian space.
\end{lemma}
For applications of this lemma see \cite{Du4}.
The existence of homogeneous
geodesics in homogeneous pseudo-Riemannian spaces (for both reductive and non reductive) was answered positively only recently by  
Z. Du\v sek in \cite{Du3}.

Two-dimensional and three-dimensional homogeneous
pseudo-Riemannian manifolds are reductive (\cite{Ca1}, \cite{Fe-Re}). Four-dimensional non reductive
homogeneous pseudo-Riemannian manifolds were classified by M.E. Fels and A.G. Renner in \cite{Fe-Re} in terms of their non reductive Lie algebras.
Their invariant pseudo-Riemannian metrics, together with their connection and
curvature, were explicitly described in by G. Calvaruso and A. Fino in \cite{Ca-Fi}.

The three-dimensional pseudo-Riemannian g.o. spaces were classified by G. Calvaruso and Marinosci in \cite{Ca-Ma2}.
In the recent work \cite{Ca-Fi-Za}, 
G. Calvaruso, A. Fino and A. Zaeim 
obtained  explicit examples of four-dimensional non reductive pseudo-Riemannian g.o. spaces.
They deduced
an explicit description in coordinates for all invariant metrics of non reductive
homogeneous pseudo-Riemannian four-manifolds.
For  those four-dimensional non reductive pseudo-Riemannia spaces which are not g.o., they determined 
 the homogeneous geodesics though a point.

\section{Two-step homogeneous geodesics}

In the  work \cite{Arv-Sou2} N.P. Souris and the author
considered a generalisation of homogeneous geodesics, namely geodesics of the form 
\begin{equation}\label{intro1}\gamma(t)=\exp(tX)\exp(tY)\cdot o, \quad X,Y\in \fr{g},\end{equation}
which we named \emph{two-step homogeneous geodesics}.  We obtained sufficient conditions on a Riemannian homogeneous space $G/K$, which  imply the existence of two-step homogeneous geodesics in $G/K$.  
A Riemannian homogeneous spaces $G/K$ such that any geodesic of $G/K$ passing through the origin is two-step homogeneous is called \emph{two-step g.o. spaces}.

Geodesics of the form (\ref{intro1}) had appeared in the work \cite{Wa} of H.C. Wang as geodesics in a semisimple Lie group $G$, equipped with a metric induced by a Cartan involution of the Lie algebra $\fr{g}$ of $G$.  Also, in \cite{Do} R. Dohira proved that if the tangent space $T_o(G/K)$ of a homogeneous space splits into submodules $\fr{m}_1,\fr{m}_2$ satisfying certain algebraic relations, and if $G/K$ is endowed with a special one parameter family of Riemannian metrics $g_c$, then all geodesics of the Riemannian space $(G/K,g_c)$ are of the form (\ref{intro1}).   
The main result of \cite{Arv-Sou2} is the following:

\begin{theorem}\label{maintheorem}\textnormal{(\cite{Arv-Sou2})}
 Let $M=G/K$ be a homogeneous space admitting a naturally reductive Riemannian metric.  Let $B$ be the corresponding inner product on $\fr{m}=T_o(G/K)$.  We assume that $\fr{m}$ admits an $\operatorname{Ad}(K)$-invariant orthogonal decomposition 

\begin{equation}\label{decco}\fr{m}=\fr{m}_1\oplus \fr{m}_2\oplus \cdots \oplus \fr{m}_s,\end{equation}

\noindent with respect to $B$.  We equip $G/K$ with a $G$-invariant Riemannian metric $g$ corresponding to the $\operatorname{Ad}(K)$-invariant positive definite inner product 
$\langle\cdot , \cdot\rangle=\lambda_1\left.B\right|_{\fr{m}_1}+\cdots +\lambda_s\left.B\right|_{\fr{m}_s}$,  $\lambda_1,\dots,\lambda_s>0$.
 If $(\fr{m}_a,\fr{m}_b)$ is a pair of submodules in the decomposition (\ref{decco}) such that 

\begin{equation}\label{geod}[\fr{m}_a,\fr{m}_b]\subset \fr{m}_a,\end{equation}

\noindent then any geodesic $\gamma$ of $(G/K,g)$ with $\gamma(0)=o$ and $\dot{\gamma}(0)\in \fr{m}_a\oplus \fr{m}_b$, is a two-step homogeneous geodesic.
In particular, if $\dot{\gamma}(0)=X_a+X_b \in \fr{m}_a\oplus \fr{m}_b$, then for every $t\in \mathbb R$ this geodesic is given by 

\begin{equation*}\gamma(t)=\exp t(X_a+\lambda X_b)\exp t(1-\lambda)X_b\cdot o,
\quad \mbox{where}\  \lambda=\lambda_b/\lambda_a.
\end{equation*}

\noindent
Moreover, if either 
$\lambda_a=\lambda_b$ or
$[\fr{m}_a,\fr{m}_b]=\left\{ {0} \right\}$ holds,
then $\gamma$ is a homogeneous geodesic, that is $\gamma(t)=\exp t(X_a+X_b)\cdot o$, for any $t\in \mathbb R$.
\end{theorem}
 
The following corollary provides a method to obtain many examples of two-step g.o. spaces.

\begin{corollary}\label{criterion}Let $M=G/K$ be a homogeneous space admitting a naturally reductive Riemannian metric.  Let $B$ be the corresponding inner product of $\fr{m}=T_o(G/K)$.  We assume that $\fr{m}$ admits an $\operatorname{Ad}(K)$-invariant, $B$-orthogonal decomposition
$\fr{m}=\fr{m}_1\oplus\fr{m}_2$, such that 
$[\fr{m}_1,\fr{m}_2]\subset \fr{m}_1$.  
Then $M$ admits an one-parameter family of $G$-invariant Riemannian metrics $g_{\lambda}$, $\lambda \in \mathbb R^+$, such that $(M,g_{\lambda})$ is a two-step g.o. space.
Each metric $g_{\lambda}$ corresponds to an $\operatorname{Ad}(K)$-invariant positive definite inner product on $\fr{m}$ of the form 

\begin{equation*}\langle \ ,\ \rangle=\left.B\right|_{\fr{m}_1}+\lambda \left.B\right|_{\fr{m}_2}, \ 
\lambda >0.\end{equation*}
\end{corollary}

The above Corollary \ref{criterion} 
is a generalisation of Dohira's result \cite{Do}.

The main tool for the proof of Theorem \ref{maintheorem} is the following proposition.

\begin{proposition} \label{generalform}\textnormal{(\cite{Arv-Sou1})}Let $M=G/K$ be a homogeneous space and $\gamma:\mathbb R\rightarrow M$ be the curve $\gamma(t)=\exp(tX)\exp(tY)\exp(tZ)\cdot o$, where $X,Y,Z\in \frak{m}$.  Let $T:\mathbb R \rightarrow \operatorname{Aut}(\fr{g})$ be the map $T(t)=\operatorname{Ad}(\exp(-tZ)\exp(-tY))$. Then $\gamma$ is a geodesic in $M$ through $o=eK$ if and only if for any $W\in \fr{m}$, the function $G_W:\mathbb R\rightarrow \mathbb R$ given by 
\begin{equation*}G_W(t)=\langle \ (TX)_{\frak{m}}+(TY)_{\frak{m}}+Z_{\frak{m}},[W, TX+TY+Z]_{\frak{m}}\ \rangle+\langle \ W,[TX, TY+Z]_{\frak{m}}+[TY, Z]_{\frak{m}}\ \rangle,\end{equation*}
is identically zero, for every $t\in \mathbb R$.\end{proposition}

The above proposition is a new tool towards the study of geodesics consisting of more than one exponential factors.  In fact, for
$X=Y=0$ we obtain
 Lemma \ref{L1} of  Kowalski and Vanhecke.

A natural application of Corollary \ref{criterion} is for total spaces of homogeneous Riemannian submersions, as shown below.

\begin{proposition}\label{su}Let $G$ be a Lie group admitting a bi-invariant Riemannian metric and let $K,H$ be closed and connected subgroups of $G$, such that $K\subset H\subset G$.  Let $B$ be the $\operatorname{Ad}$-invariant positive definite inner product on the Lie algebra $\fr{g}$ corresponding to the bi-invariant metric of $G$.  We identify each of the spaces $T_o(G/K),T_o(G/H)$ and $T_o(H/K)$ with corresponding subspaces $\fr{m},\fr{m}_1$ and $\fr{m}_2$ of $\fr{g}$, such that $\fr{m}=\fr{m}_1\oplus \fr{m}_2$.  We endow $G/K$ with the $G$-invariant Riemannian metric $g_{\lambda}$ corresponding to the $\operatorname{Ad}(K)$-invariant positive definite inner product  
$\langle\cdot \ ,\ \cdot\rangle=\left.B\right|_{\fr{m}_1}+\lambda \left.B\right|_{\fr{m}_2}, \quad \lambda>0$.
 Then $(G/K,g_{\lambda})$ is a two-step g.o. space.

\end{proposition}

\begin{example}\textnormal(\cite{Arv-Sou1}) The odd dimensional sphere $\mathbb S^{2n+1}$ can be considered as the total space of the homogeneous Hopf bundle 
$\mathbb S^1 \rightarrow \mathbb S^{2n+1} \rightarrow \mathbb CP^n$.
Let $g_1$ be the standard metric of $\mathbb S^{2n+1}$.  We equip $\mathbb S^{2n+1}$ with an one parameter family of metrics $g_{\lambda}$, which ``deform" the standard metric along the Hopf circles $\mathbb S^1$.  
By setting $G=\U(n+1)$, $K=\U(n)$ and $H=\U(n)\times \U(1)$, the Hopf bundle corresponds to the fibration
$H/K \rightarrow G/K \rightarrow G/H$.

Since $\U(n+1)$ is compact, it admits a bi-invariant metric corresponding to an $\operatorname{Ad}(\U(n+1))$-invariant positive definite inner product $B$ on $\fr{u}(n+1)$.  We identify each of the spaces $T_o\mathbb S^{2n+1}= T_o(G/K),T_o\mathbb CP^n= T_o(G/H)$, and $T_o\mathbb S^1= T_o(H/K)$ with corresponding subspaces $\fr{m},\fr{m}_1$, and $\fr{m}_2$ of $\fr{u}(n+1)$.  The desired one parameter family of metrics $g_{\lambda}$ corresponds to the one parameter family of positive definite inner products 
$\langle \ ,\ \rangle=\left.B\right|_{\fr{m}_1}+\lambda\left.B\right|_{\fr{m}_2}$, $\lambda>0$
\noindent on $\fr{m}=\fr{m}_1\oplus \fr{m}_2$.  

  Then Proposition \ref{su} implies that $(\mathbb S^{2n+1},g_{\lambda})$ is a two-step g.o. space.  In particular, let $X\in T_o\mathbb S^{2n+1}$.  Then the unique geodesic $\gamma$ of  $(\mathbb S^{2n+1},g_{\lambda})$ with $\gamma(0)=o$ and $\dot{\gamma}(0)=X$, is given by
$\gamma(t)=\exp t(X_1+\lambda X_2) \exp t(1-\lambda)X_2\cdot o$,
where $X_1$, $X_2$ are the projections of $X$ on $\fr{m}_1= T_o\mathbb CP^n$ and $\fr{m}_2= T_o\mathbb S^1$ respectively.\\
Note that if $\lambda=1+\epsilon$, $\epsilon >0$, then the metrics $g_{1+\epsilon}$ are Cheeger deformations of the natural metric $g_1$.  
\end{example}

By using Proposition \ref{criterion} it is possible to construct various classes of two-step g.o spaces.  The recipe is the following:

\noindent
(i) Let $G/K$ be  a homogeneous space with reductive decomposition $\fr{g}=\fr{k}\oplus\fr{m}$ admitting a naturally reductive metric corresponding to a positive definite inner product $B$ on $\fr{m}$.

\noindent
  (ii) We consider an $\operatorname{Ad}(K)$-invariant, orthogonal decomposition 
$\fr{m}=\fr{n}_1\oplus \cdots \oplus \fr{n}_s$ with respect to $B$. 

\noindent
(iii) We separate the submodules $\fr{n}_i$ into two groups as
$\fr{m}_1=\fr{n}_{i_1}\oplus \cdots \oplus \fr{n}_{i_n} \quad \makebox{and} \quad \fr{m}_2=\fr{n}_{i_{n+1}}\oplus \cdots \oplus \fr{n}_{i_s}$,
so that $[\fr{m}_1,\fr{m}_2]\subset \fr{m}_1$.
The decomposition $\fr{m}=\fr{m}_1\oplus \fr{m}_2$ is $\operatorname{Ad}(K)$-invariant and orthogonal with respect to $B$.  

\noindent
(iv) Then Corollary \ref{criterion}  implies that $G/K$ admits an one parameter family of metrics $g_{\lambda}$ so that $(G/K,g_{\lambda})$ is a two-step g.o. space.

In \cite{Arv-Sou2} we
 applied the above recipe to the following classes of homogeneous spaces:

\noindent
1)\quad Lie groups with bi-invariant metrics, equipped with an one parameter family of left-invariant metrics.

\noindent
2)\quad Flag manifolds equipped with certain one parameter families of diagonal metrics.

\noindent
3)\quad Generalized Wallach spaces equipped with three different types of diagonal metrics (thus recovering some results from \cite{Arv-Sou1}.

\noindent
4)\quad $k$-symmetric spaces where $k$ is even, endowed with an one parameter family of diagonal metrics.

\section{Some open problems}
It seems that the target for a complete classification of homogeneous g.o. spaces in any dimension greater than seven is far for being accomplished. In dimension seven there are several examples but a complete classification is still unknown.
However, as shown in the present paper, for some large classes of homogeneous spaces it is possible to obtain some necessary conditions for the g.o. property.  These conditions are normally imposed by  the
 special Lie theoretic structure of corresponding homogeneous space.
Also, the problem of an explicit description of homogeneous geodesics for spaces which are not g.o., is not trivial either.
Eventhough it is mathematically simple, it requires high computational complexity.

Another difficulty that one faces, is 
 to show that the g.o. property of a homogeneous space 
$(M=G/K,g)$ does not depend on the representation as a coset space and on the $\Ad(K)$-invariant decomposition  $\g =\kk \oplus\mm $, that is the space is in no way naturally reductive. 
Therefore, we often stress that we study $G$-g.o. spaces.

Also, it would be interesting to see how various results about Riemannian manifolds could be adjusted to pseudo-Riemannian manifolds, such as Propositions \ref{genlemma}, \ref{generalform}.

Concerning generalizations of the g.o. propery, we have introduced the concept of a two-step homogeneous geodesic and two-step g.o. space.
We conjecture that a search for  three-step (or more) homogeneous geodesics   reduces to two-step homogeneous geodesics.

\end{document}